%&amstex
\input amstex
\input amsppt.sty
\magnification=\magstep1
\hsize=30truecc
\vsize=22.2truecm
\baselineskip=16truept
\NoBlackBoxes
\nologo \pageno=1
\TagsOnRight
\topmatter
\def\Z{\Bbb Z}

\def\l{\left}
\def\r{\right}
\def\bg{\bigg}
\def\({\bg(}
\def\[{\bg[}
\def\){\bg)}
\def\]{\bg]}
\def\t{\text}
\def\f{\frac}
\def\mo{\roman{mod}}

\def\em{\emptyset}
\def\se {\subseteq}

\def\sm{\setminus}

\def\eq{\equiv}

\def\ls{\leqslant}
\def\gs{\geqslant}

\topmatter \hbox{Adv. Appl. Math. 35(2005), no.\,2, 182--187.}
\bigskip
\title {On odd covering systems with distinct moduli}\endtitle
\author {Song Guo and Zhi-Wei Sun*}\endauthor
\leftheadtext{Song Guo and Zhi-Wei Sun}
\affil Department of Mathematics, Nanjing University
\\Nanjing 210093, People's Republic of China
\\ {\tt guosong77\@sohu.com}\ \ \quad \ \ \quad\tt{zwsun\@nju.edu.cn}
\endaffil
\date Received 10 December 2004; accepted 22 January 2005\enddate
\abstract A famous unsolved conjecture of P. Erd\H os and J. L.
Selfridge states that there does not exist a covering system
 $\{a_s(\mo\ n_s)\}_{s=1}^{k}$ with the moduli $n_1,\ldots,n_k$ odd,
 distinct and greater than one. In this paper we show that if such a covering
system $\{a_s(\mo\ n_s)\}_{s=1}^{k}$ exists with
$n_1,\ldots,n_k$ all square-free, then the least common multiple
of $n_1,\ldots,n_k$ has at least 22 prime divisors.
\endabstract
\thanks  2000 {\it Mathematics Subject Classification}.
Primary 11B25; Secondary 11A07.
\newline\indent *This author is responsible for all communications,
and supported by the National Science Fund for Distinguished Young
Scholars (No. 10425103) and the Key Program of NSF (No. 10331020)
in P. R. China.
\endthanks
\endtopmatter
\document
\hsize=30truecc
\vsize=22.2truecm
\baselineskip=16truept

\heading 1. Introduction\endheading
 For $a\in \Z$ and $n\in\{1,2,3,\ldots\}$, we simply let $a(n)$ denote the residue class
 $$a(\mo\ n)=\{a+nx:\, x\in\Z\}.$$
In the early 1930s P. Erd\H os called a finite system
$$A=\{a_s(n_s)\}_{s=1}^k\tag $*$ $$
of residue classes a {\it covering system} if
$\bigcup_{s=1}^ka_s(n_s)=\Z$. Clearly $(*)$ is a covering system
if and only if it covers $0,1,\ldots,N_A-1$ where
$N_A=[n_1,\ldots,n_k]$ is the least common multiple of the moduli
$n_1,\ldots,n_k$.

 Here are two covering systems with distinct moduli
constructed by Erd\H os:
$$\{0(2),\ 0(3),\ 1(4),\ 5(6),\ 7(12)\},$$
$$\align \{0(2),\ &0(3),\ 0(5),\ 1(6),\ 0(7),\ 1(10),\ 1(14),\ 2(15),
           \\&2(21),\ 23(30),\ 4(35),\ 5(42),\ 59(70),\ 104(105)\}.
\endalign$$

 Covering systems have been investigated by various number theorists and
combinatorists, and many surprising applications have been found.
(See [3], [4] and [7].)

A covering system with odd moduli
is said to be an {\it odd covering system}.
Here is a well-known open problem in the field (cf. [3]).

\proclaim{Erd\H os-Selfridge Conjecture} There does not exist an
odd covering system with the moduli distinct and greater than one.
\endproclaim

 In 1986-1987, by a lattice-geometric method,
 M. A. Berger, A. Felzenbaum and A. S. Fraenkel ([1] and [2]) obtained
some necessary conditions for system $(*)$ to be an odd covering
system with $1<n_1<\cdots<n_k$, one of which is the inequality
$$\prod_{t=1}^{r}\frac{p_t-1}{p_t-2}-\sum_{t=1}^{r}\frac{1}{p_t-2}>2,$$
where $p_1,\ldots,p_r$ are the distinct prime divisors of $N_A$.
They also showed that if $(*)$ is an odd covering system with
$n_1,\ldots,n_k$ square-free, distinct and greater than one, then
the above inequality can be improved as follows:
$$\prod_{t=1}^{r}\frac{p_t}{p_t-1}-\sum_{t=1}^{r}\frac{1}{p_t-1}\gs2$$
and consequently $r\gs 11$. This was also deduced
by the second author [6] in a simple way.

 In 1991, by a complicated sieve method, R. J. Simpson
 and D. Zeilberger [5] proved that
if $(*)$ is an odd covering system with $n_1,\ldots,n_k$
square-free, distinct and greater than one, then $N_A$ has at
least 18 prime divisors.

  In this paper we obtain further improvement in this direction
by a direct argument.

 \proclaim {Theorem 1} Suppose that $(*)$ is an odd covering system
 with $1<n_1<\cdots<n_k$.
 If $N_A=[n_1,\ldots,n_k]$ is square-free,
then it has at least $22$ prime divisors.
\endproclaim

 In contrast with the Erd\H os-Selfridge conjecture,
 recently the second author [8] showed that
 if $(*)$ is a covering system with $1<n_1<\cdots<n_k$
 then it cannot cover every integer an odd number of times.

\heading {2. Proof of Theorem 1}\endheading

For convenience we let $[a,b]=\{x\in\Z:\, a\ls x\ls b\}$ for any $a,b\in\Z$.

Assume that $N=N_A=p_1\cdots p_r$ where $p_1<\cdots<p_r$
are distinct odd primes.
For each $t\in[1,r]$, we set
$$d_t=\l\lfloor\f35(t-1)\r\rfloor$$
(where $\lfloor\cdot\rfloor$
is the greatest integer function), and define
$$M_t=\cases\{p_ip_t:\,1\ls i\ls d_t\}&\t{if}\ t\ls 8,
\\\{p_ip_t:\,1\ls i\ls d_t\}\cup\{p_1p_2p_t,p_1p_3p_t\}&\t{if}\ t\gs 9.
\endcases$$
Note that $d_1=d_2=0$ and hence $M_1=M_2=\em$.

For $s\in[1,k]$ let
$$n_s'= \cases
p_t&\t{if}\ n_s\in M_t\ \t{for some}\ t,\\
n_s&\t{otherwise}.
\endcases$$
Since $n_s'\mid n_s$, we have $a_s(n_s)\se a_s(n_s')$. Thus
$A'$=$\lbrace a_s(n_s')\rbrace _{s=1}^{k}$ is also an odd covering
system. Let
$$\bar I = [1,k]\backslash \bigcup_{t=1}^{r}I_t
\ \ \ \ \t{where}\ I_t=\{1\ls s\ls k:\, n_s'=p_t\}.$$ Then
$$\bigcup_{s \in \bar I}a_s(n_s')\supseteq [0,N-1]
\backslash \bigcup_{t=1}^{r}\bigcup_{s \in I_t}a_s(n_s')
=\bigcap_{t=1}^{r}\bg([0,N-1] \backslash  \bigcup_{s\in
I_t} a_s(n_s')\bg).$$

For each $t\in[1,r]$, clearly $|I_t|\ls d_t+1<p_t$ if $t\ls 8$,
and $|I_t|\ls d_t+3$ otherwise. Observe that $d_t\ls 3(p_t-1)/5
<p_t-3$ if $t\gs 9$. So there is a subset $R_t$ of $[0,p_t-1]$
satisfying the following conditions:

(a) $|R_t|=p_t-1-d_t$ if $t\ls 8$, and $|R_t|=p_t-3-d_t$ if $t\gs 9$;

(b) $x\not\equiv a_s(\mo\ p_t)$ for any $x\in R_t$ and $s\in I_t$.

 Define
$$X=\{x\in[0,N-1]:\ \t{the remainder of}\ x\ \mo\ p_t
\ \t{lies in}\ R_t\ \t{for}\ t\in[1,r]\}.$$
Then $|X|=\prod_{t=1}^r|R_t|$ by the Chinese Remainder Theorem, also
$$X\subseteq
\bigcap_{t=1}^{r}\bg([0,N-1] \backslash \bigcup_{s\in I_t}
a_s(n_s')\bg) \subseteq \bigcup_{s\in \bar I}a_s(n_s')
=\bigcup_{s\in \bar I}a_s(n_s)$$ and hence
$X=\bigcup_{s\in J}X_s$, where
$$X_s=X\cap a_s(n_s)\ \ \t{and}\ \ J=\{s\in \bar I:\, X_s\not=\em\}.$$
For each $s\in J$, the set $X_s$ consists of
those $x\in[0,N-1]$ for which $x\eq a_s\ (\mo\ p_t)$ if $p_t\mid
n_s$, and $x\eq r_t\ (\mo\ p_t)$ for some $r_t\in R_t$ if
$p_t\nmid n_s$. Thus, by the Chinese Remainder Theorem,
$$|X_s|=\prod\Sb 1\ls t\ls r\\p_t\nmid n_s\endSb|R_t|
=|X|\prod\Sb 1\ls t\ls r\\p_t\mid n_s\endSb|R_t|^{-1} \ \ \t{for
all}\ s\in J.$$

  Let $a_0\in X$, $n_0=p_1p_2$ and $X_0=X\cap a_0(n_0)$.
 Again by the Chinese Remainder Theorem,
 $$|X_0|=\prod_{2<t\ls r}|R_t|=|X|\prod\Sb 1\ls t\ls r\\p_t\mid n_0\endSb
 |R_t|^{-1}.$$
 Let $j=0$ if $n_0\not\in\{n_s:\,s\in J\}$, and let $j$
 be the unique element of $J$ with $n_j=n_0$ if $n_0\in\{n_s:\,s\in J\}$.
Set $J_0=\{s\in J:\,(n_s,n_0)=1\}$. Then
$$\aligned|X|=\bg|\bigcup_{s \in J\cup\{j\}}X_s\bg|
\ls&\sum_{s \in J\sm (J_0\cup\{j\})}|X_s|+\bg|X_j
\cup \bigcup_{s\in J_0}X_s\bg|
\\\ls&\sum_{s\in J\sm (J_0\cup\{j\})}|X_s|+|X_j|
+\sum_{s\in J_0}|X_s\backslash X_j|
\\=&\sum_{s \in J\sm (J_0\cup\{j\})}|X_s|+|X_j|
+\sum_{s\in J_0}(|X_s|-|X_s\cap X_j|)
\endaligned$$
and so
$$|X|\ls \sum_{s \in J\cup\{j\}}|X_s|-\sum_{s\in J_0}|X_s\cap X_j|.$$
If $s\in J_0$, then
$X_s\cap X_j$ consists of those $x\in[0,N-1]$ for which
$x\eq a_j\ (\mo\ n_j)$, $x\eq a_s\ (\mo\ n_s)$, and
$x\eq r_t\ (\mo\ p_t)$ for some $r_t\in R_t$ if $p_t\nmid n_jn_s$,
therefore
$$|X_s\cap X_j|=\prod\Sb 1\ls t\ls r\\p_t\nmid n_jn_s\endSb
|R_t|=|X|\prod\Sb 1\ls t\ls r\\p_t\mid n_0n_s\endSb |R_t|^{-1}.$$

Set
$$D_1=\{d>1:\,d\mid N\}\sm\bg(\{p_1,\ldots,p_r\}\cup\bigcup_{t=1}^r M_t\bg),$$
and $$D_2=\{n_s:\, s\in J\cup\{j\}\} \ \t{and}\
D_3=\{d\in D_1:\,(d,n_0)=1\}.$$
If $s\in J$, then $n_s'\not=p_t$ for any $t\in[1,r]$,
and thus $n_s=n_s'\in D_1$. Since $d_2=0$, we also have $n_j=p_1p_2\in D_1$.
Therefore $D_2\se D_1$, and so $D_2 \cap D_3$ coincides with
$D_4=\{n_s:\,s\in J_0\}$.

Let
$$x_t=|R_t|^{-1}\ls1\quad\t{for}\ t=1,\ldots,r,$$
and
$$I(d)=\{1\ls t\ls r:\, p_t\mid d\}$$
for any positive divisor $d$ of $N$. Observe that
$$\align&\sum_{d\in D_1\sm D_2}\prod_{t\in I(d)}x_t
-x_1x_2\sum_{d\in D_3\sm D_4}\prod_{t\in I(d)}x_t
\\=&\sum_{d\in D_1\sm(D_2\cup D_3)}\prod_{t\in I(d)}x_t
+(1-x_1x_2)\sum_{d\in D_3\sm D_4}\prod_{t\in I(d)}x_t\gs0.
\endalign$$
Thus
$$\align|X|
\ls&\sum_{s \in J\cup\{j\}}|X_s|-\sum_{s\in J_0}|X_s\cap X_j|
\\=&\sum_{d\in D_2}|X|\prod_{t\in I(d)}x_t
-\sum_{d \in D_4}|X|x_1x_2\prod_{t\in I(d)}x_t
\\\ls &|X|\(\sum_{d\in D_1} \prod_{t\in I(d)}x_t
-x_1x_2\sum_{d \in D_3}\prod_{t\in I(d)}x_t\).
\endalign$$

 Since $d_1=d_2=0$ and $d_t<3$ for $t<6$, by the above we have
$$\aligned1\ls&\sum\Sb I\se [1,r]\\|I|>1\endSb\prod_{t \in I}x_t
-\sum_{t=1}^{r}\sum_{1\ls i\ls d_t}x_ix_t
-\sum_{9\ls t\ls r}(x_1x_2x_t+x_1x_3x_t)
\\&-x_1x_2\(\sum\Sb I\se [3,r]\\|I|>1\endSb\prod_{t
\in I}x_t -\sum_{3\ls t\ls r}\sum_{3\ls i\ls d_t}x_ix_t\)
\\=&\prod_{t=1}^{r}(1+x_t)-1-\sum_{t=1}^{r}x_t
-\sum_{t=3}^{r}\sum_{i=1}^{d_t}x_ix_t
-\sum_{9\ls t\ls r}(x_1x_2x_t+x_1x_3x_t)
\\&-x_1x_2\(\prod_{t=3}^{r}(1+x_t)-1-\sum_{t=3}^{r}x_t
-\sum_{t=6}^{r}\sum_{i=3}^{d_t}x_ix_t\).
\endaligned$$
It follows that $f(x_1,\ldots,x_r)\gs 2$,
where
$$\aligned f(x_1,...,x_r)=&(1+x_1+x_2)\prod_{t=3}^{r}(1+x_t)
-\sum_{t=1}^r x_t+x_1x_2-\sum_{t=3}^{r}\sum_{i=1}^{d_t}x_ix_t
\\&+x_1x_2\sum_{t=3}^{8}x_t-x_1x_3\sum_{9\ls t\ls r}x_t
+x_1x_2\sum_{t=6}^{r}\sum_{i=3}^{d_t}x_ix_t
\endaligned$$
can be written in the form
$\sum_{i_1,\ldots,i_r}c_{i_1,\ldots,i_r}x_1^{i_1}\cdots x_r^{i_r}$
with $c_{i_1,\ldots,i_r}\gs0$.

Let $q_1=3<\cdots<q_r$ be the first $r$ odd primes.
For each $t\in[1,r]$, as $p_t\gs q_t$ we have $x_t\ls x_t'$, where
$$x_t'=\cases(q_t-d_t-1)^{-1}&\t{if}\ 1\ls t\ls 8,
\\(q_t-d_t-3)^{-1}&\t{if}\ 9\ls t\ls r.\endcases$$
Thus
$$f(x_1',\ldots,x_r')\gs f(x_1,\ldots,x_r)\gs2.$$
By computation through computer we find that
$$f(x_1',\ldots,x_{21}')=1.995\cdots<2,$$
therefore $r\not=21$. (This is why we define $d_t$ and $M_t$
in a somewhat curious way.)

In the case $r<21$, we let $p_{r+1}<\cdots<p_{21}$ be distinct primes
greater than $p_r$, and then
$$\Cal A=\{a_1(n_1),\ldots,a_k(n_k),0(p_{r+1}),\ldots,0(p_{21})\}$$
forms an odd covering system with $N_{\Cal A}$ square-free and
having exactly 21 distinct prime divisors. This is impossible by the above.

 Now we can conclude that $r\gs 22$ and this completes the proof.

\enddocument